\documentclass[11pt,oneside,english,reqno]{amsart}
\usepackage[T1]{fontenc}
\usepackage[utf8]{inputenc}
\usepackage[a4paper]{geometry}
\usepackage{mathrsfs}
\usepackage{amstext}
\usepackage{amsthm}
\usepackage{amssymb}
\usepackage{color}
\usepackage{hyperref}
\usepackage{enumerate}
\usepackage{verbatim} 
\usepackage{stmaryrd}
\usepackage{enumitem}
\usepackage{graphicx}

\makeatletter
\numberwithin{equation}{section}
\numberwithin{figure}{section}
\theoremstyle{plain}

\theoremstyle{definition}

\theoremstyle{remark}

\theoremstyle{plain}

\theoremstyle{plain}

\theoremstyle{plain}

\makeatother

\usepackage{babel}
\providecommand{\corollaryname}{Corollary}
\providecommand{\definitionname}{Definition}
\providecommand{\lemmaname}{Lemma}
\providecommand{\propositionname}{Proposition}
\providecommand{\remarkname}{Remark}
\providecommand{\theoremname}{Theorem}

\begin{document}
   
\title{Can probability theory really help tame problems in mathematical hydrodynamics?}
\author{Martina Hofmanov\'a, Florian Bechtold}
\email{hofmanova@math.uni-bielefeld.de, fbechtold@math.uni-bielefeld.de}
\maketitle
\begin{abstract}
Recent years have seen spectacular progress in the mathematical study of hydrodynamic equations. Novel tools from convex integration in particular prove extremely versatile in establishing non-uniqueness results. Motivated by this 'pathological' behavior of solutions in the deterministic setting, stochastic models of fluid dynamics have enjoyed growing interest from the mathematical community. Inspired by the theory of 'regularization by noise', it is hoped for that stochasticity might help avoid 'pathologies' such as non-uniqueness of weak solutions. Current research however shows that convex integration methods can prevail even in spite of random perturbations.  

\end{abstract}

\section{Some recent developments in mathematical fluid dynamics}
The main object of interest in mathematical fluid dynamics consists in the Navier-Stokes equations
\begin{equation}
\begin{split}
    \label{NS}
    \frac{\partial}{\partial t}u(t,x)+(u(t,x)\cdot \nabla )u(t,x)+\nabla p(t,x)&=\nu \Delta u(t,x)\\
    \mbox{div}\ u(t,x)&=0\\
    u(0, x)&=u_0(x),
\end{split}
\end{equation}
which model the evolution of the velocity field $u:[0,T]\times\mathbb{R}^3\to \mathbb{R}^3$ and pressure $p:[0,T]\times\mathbb{R}^3\to \mathbb{R}$ of an incompressible fluid of viscosity $\nu>0$. In the case of ideal fluids ($\nu=0$) the equations \eqref{NS} are called Euler equations. While ubiquitous in applications  e.g. in engineering, a complete mathematical understanding of the Navier-Stokes equations remains elusive. Arguably the most famous problem consists in the study of global smooth solutions to \eqref{NS} in three space dimensions. 

While this Millennium problem remains open, other long standing questions of well-posedness were successfully answered in the past years: Buckmaster and Vicol \cite{Buckmaster2019b} established non-uniqueness of weak solutions with finite kinetic energy. More precisely, the authors showed that for any prescribed smooth and non-negative function $e$ there is a weak solution whose kinetic energy is given by $e$: Translated into physics, such solutions would for instance model the case of a glass of water at rest which spontaneously begins to stir itself. As such a behavior clearly contradicts basic physical principles, such solutions are sometimes called 'pathological': while they are proper weak solutions to \eqref{NS} in the mathematical sense, they do not actually model the physical behavior of a fluid.

This particular and many related results rely on the method of convex integration introduced to fluid dynamics by De Lellis and Sz\'ekelyhidi \cite{DeLellis2009, deLellis2008, DeLellis2012}. Before its application to the Navier-Stokes equations, this method has already led to a number of groundbreaking results concerning the incompressible Euler equations, culminating in the proof of Onsager’s conjecture by Isett \cite{Isett2018} and by Buckmaster, De Lellis, Sz\'ekelyhidi and Vicol \cite{Buckmaster2018onsager}. We refer to the excellent review articles by Buckmaster and Vicol \cite{Buckmaster2020, Buckmaster2020b} for a gentle introduction to convex integration and further references. 

From the above, one infers that weak solutions are too large of a class of solutions to consider for the Navier-Stokes equations, as this class contained the above mentioned 'pathological' examples. Yet, this does not sort out the possibility of uniqueness in a more restrictive class of solutions subject to further constraints. A natural way to narrow the class of solutions would consist in imposing further physical constraints, such as energy inequalities. Very recently however, Albritton, Bru{\'e} and Colombo  \cite{colomboleray} were able to extend the above non-uniqueness results even to such Leray solutions, i.e. weak solutions satisfying an energy inequality, for forced Navier-Stokes equations. \\
\\
In light of these exciting but rather 'negative' developments, the question of a possible regularizing
effect provided by a suitable stochastic perturbation becomes even more prominent. From a physical perspective, stochastic perturbations have been considered for example to account for microscopic thermal fluctuations \cite{landau2013fluid} or as a model for turbulence \cite{kraichnan}. From a mathematical perspective, it then becomes natural to ask if such stochastic perturbations give rise to a better behaved system, e.g.: Can the stochastic problem help avoid 'pathologies' apparent in the deterministic system as established by convex integration methods and thus restore uniqueness of weak/Leray solutions? 

\section{A glimps at regularization by noise}

A natural first question one could ask at this point is: Why should the addition of some stochasticity actually improve things? To address this question, let us discuss the heuristics of a simple concrete example. As is well known ordinary differential equations of the form 
\begin{equation}
    \frac{d}{dt}y=f(y), \qquad y(0)=y_0
    \label{ode}
\end{equation}
admit a unique global solution, provided $f$ is Lipschitz continuous and of linear growth. If $f$ is less regular, assuming for example $f$ to be only $\alpha$-H\"older continuous for some $\alpha\in (0,1)$ uniqueness might fail. The most canonical stochastic version of \eqref{ode} consists in
\begin{equation}
    \frac{d}{dt}y=f(y)+\frac{d}{dt}w, \qquad y(0)=y_0
    \label{sde}
\end{equation}
where $w$ is some stochastic process modeling random perturbations of \eqref{ode}. In the case of $w$ being a Brownian motion\footnote{This means that $\frac{d}{dt}w$ is a white noise i.e. arguably the most canonical random perturbation conceivable.} it was shown already in the seminal works of Zvonkin \cite{Zvonkin_1974} and Veretennikov \cite{Veretennikov_1981} that bounded and measurable drifts $f$ are sufficient to ensure well posedness, i.e. existence and uniqueness of solutions to \eqref{sde}. These results were later extended to drifts satisfying a quantified integrability condition by Krylov and R\"ockner \cite{Krylov2004}. An excellent review on regularization by Brownian motion  can be found in \cite{Flandoli2011}.

Starting with the seminal work by Gubinelli and Catellier \cite{gubicat} the particular approach of 'pathwise regularization by noise' has enjoyed considerable attention in recent years \cite{galeati2020noiseless}, \cite{harang2020cinfinity}, \cite{galeati2020regularizationn}, \cite{bechtold}. This approach that typically investigates regularization due to very rough noises $w$ instead of Brownian motion is based on the simple observation that upon the transformation $x:=y-w$ and integration, \eqref{sde} reads
\begin{equation}
    x_t=x_0+\int_0^tf(x_s+w_s)ds.
    \label{integrated eqn}
\end{equation}
If one assumes the oscillations of the noise $w$ to dominate the oscillations of $x$ one would expect that on small time scales $t-s\ll1$, the above integral should behave as 
\begin{equation}
    \int_s^tf(x_r+w_r)dr\simeq \int_s^t f(x_s+w_r)dr.
    \label{approximation}
\end{equation}
As $w$ is highly oscillatory, one heuristically expects the integral on the right hand side to 'average out over singularities of $f$', meaning that the function 
\[
x\to \int_s^tf(x+w_r)dr
\]
is expected to enjoy better regularity properties than the function $f$. By 'gluing together' these locally regularized approximations \eqref{approximation} one can recover the global integral and thus establish a well-posedness theory for \eqref{integrated eqn} in settings of considerably less regular $f$ than required in the deterministic setting. For a nice overview of contemporary developments in this direction, we refer to \cite{luciophd}.

\section{Convex integration in the presence of noise}
Given the non-uniqueness results obtained by convex integration methods for the deterministic Navier-Stokes equations and the phenomenon of regularization by noise, one might hope that stochasticity would improve the behaviour of \eqref{NS}. 

Indeed, some positive results have
been achieved. Flandoli and Luo \cite{Flandoli2021} showed that a noise of transport type prevents a vorticity
blow-up in the Navier-Stokes equations (see also \cite{Flandoli2020, Flandoli20211} for related results). Flandoli, Hofmanov\'a, Luo and Nilssen \cite{flandohof} then
showed that the regularization is even provided by deterministic vector fields. A linear multiplicative
noise prevents the blow up of the velocity with high probability for the three dimensional Euler
and Navier–Stokes system as well, as shown by Glatt-Holtz and Vicol \cite{glattvicol} and R\"ockner, Zhu
and Zhu \cite{ROCKNER20141974}, respectively. Noise also has a beneficial impact when it comes to long time
behavior and ergodicity. Da Prato and Debussche \cite{DAPRATO2003877} obtained a unique ergodicity for three
dimensional stochastic Navier-Stokes equations with non-degenerate additive noise. The theory
of Markov selections by Flandoli and Romito \cite{Flandoli2007} provides an alternative approach which also
allowed to prove ergodicity for every Markov solution, see Romito \cite{Romito2008}.\\
\\
However, with respect to well-posedness i.e. in particular uniqueness of weak solutions to stochastic versions of \eqref{NS}, it appears that convex integration arguments prevail over several types of noises, as established by Hofmanov\'a, Zhu and Zhu in a series of works. Along the way, convex integration also serves as a remarkably useful tool to address problems in the field of stochastic partial differential equations (SPDEs):  In a first work \cite{hofmanovaconvex1} non-uniqueness in law of weak solutions to \eqref{NS} driven by additive, linear multiplicative and nonlinear noise of cylindrical type was established. In a second work \cite{Hofmanov2021convex} similar results as well as existence and non-uniqueness of strong Markov solutions are obtained for the Euler equations. A third work \cite{hofmanovaconvex2} solves the long standing problem of constructing probabilistically strong global solutions to the Navier-Stokes equation perturbed by trace-class noise, along with non-uniqueness statements for such solutions. Combining convex integration techniques with tools from paracontrolled calculus \cite{gubinelli_imkeller_perkowski_2015} for singular SPDEs, \cite{hofmanovaconvex3} is able to study the Navier-Stokes equation perturbed by space-time white noise yielding global existence and non-uniqueness of weak solutions. Remarkably, the interplay of techniques allows to push the solution theory even further and into the regimes of so called supercritical equations \cite{hofmanovaconvex4} that are inaccessible by standard theories for singular SPDEs such as paracontrolled calculus or regularity structures \cite{Hairer2014}. Finally, in the most recent work \cite{hofmanovaconvex5} a full bona fide stochastic convex integration theory is developed, leading to non-uniqueness of stationary ergodic solutions to stochastic perturbations of the Navier-Stokes and Euler equations. 

\section*{Conclusion}

The advent of convex integration techniques has revolutionized the field of fluid dynamics. New 'non-uniqueness' results stress that we are still far from a comprehensive understanding of the Navier-Stokes and Euler equations. While the addition of noise has proven to be mathematically beneficial in certain regards, convex integration techniques typically appear to prevail in the face of stochastic perturbations while  opening a new potential angle in the study of SPDEs. 

\subsection*{Acknowledgement}
The authors acknowledge support from the European Research Council (ERC) under
the European Union’s Horizon 2020 research and innovation programme (grant agreement
No. 949981).
\bibliography{main}
\bibliographystyle{amsplain}

\end{document}